 \documentclass{article}
 \usepackage{amsmath,amsthm,amssymb,mathrsfs,latexsym,amsfonts,a4wide}

\newcommand{\ts}{\otimes} 

\newcommand{\pp}{{\bf P}}

\newcommand{\Z}{{\bf Z}}

\newcommand{\oo}{{\cal O}}
\newcommand{\ii}{{\cal I}}

\newcommand{\ra}{\rightarrow}

\newcommand{\mult}{{\rm mult}}

\newtheorem{theorem}{Theorem}[section]
\newtheorem{lemma}[theorem]{Lemma}

\newtheorem{cor}[theorem]{Corollary}

\newtheorem{definition}[theorem]{Definition}

\numberwithin{equation}{section}
\newcommand{\rn}{ }

\begin{document}
\title{Seshadri Constants and Interpolation on Commutative Algebraic Groups}
\author{St\'ephane Fischler and Michael Nakamaye}
\maketitle

\newcommand{\pb}[1]{{\bf #1}}
\medskip
\medskip

\begin{abstract}
In this article we study interpolation estimates on a special class of compactifications
of commutative algebraic groups constructed by Serre.  We obtain a large quantitative
improvement over previous results due to Masser and the first author and our main result
has the same level of accuracy as the best known multiplicity estimates.  The improvements
come both from using special properties of the compactifications which we
consider and from a different approach based upon Seshadri constants and vanishing
theorems. 
\end{abstract}

\medskip

\section*{Introduction}

Suppose $X$ is a smooth projective variety, $L$ a line bundle on $X$, and
$P_1,\ldots, P_n$ points in $X$.  The basic interpolation  question
asks to determine, with respect to the positivity of $L$ and the
geometry of the collection of points $\{P_i\}$, whether or not the natural
evaluation map
$$
H^0(X,L) \ra H^0\left(X, L \ts \oplus_{i=1}^n \oo_X/m_{P_i}^r\right)
$$
is surjective for a given positive integer $r$,  where $m_x \subset \oo_X$ 
is the maximal ideal sheaf corresponding to the point $x \in X$.
More generally, the ideals $m_{P_i}^r$ may be replaced with
ideals which only measure vanishing in ``certain directions''.  
A first result of this type, in the context of commutative algebraic groups and with $r =1$,  was established by Masser
\cite{M} and generalized by  the first author \cite{Fi} to allow an arbitrary
value of $r$ and also  substantial flexibility in the 
choice of the ideal sheaves that can be used instead of $m_{P_i}^r $;  in this paper we further refine these estimates.  The importance
of this question, in addition to being a fundamental issue of positivity 
in algebraic geometry, lies in its potential applications in transcendence
theory. Such a result can replace a zero (or multiplicity) estimate in a transcendence proof, 
by constructing an auxiliary functional instead of an auxiliary function 
(see \cite{W1},  \cite{W2}, \cite{W3}, and \cite{W4}).  From a Diophantine point of view, this new strategy may lead
to new results in transcendence or Diophantine approximation because the parameters are chosen in a different way.
In addition, our new proof of an interpolation estimate is geometric in nature, using
Seshadri exceptional subvarieties; in particular, we are able, as in \cite{NR}, to relate  the algebraic
subgroup responsible for a possible defect in this estimate to a Seshadri exceptional subvariety.

\medskip

In order to state our results, we require a definition of Seshadri constants.  These
numerical invariants attached to a line bundle and a finite set of points measure
the answer to the interpolation question as we will see below. Greater details will be found in \S \ref{secsesh}.
\begin{definition}
Suppose $X$ is a smooth projective variety, $L$ an ample line 
bundle on $X$, and $\Omega \subset X$ a finite set.
We define the Seshadri constant of $L$ at $\Omega$ by
$$
\epsilon(\Omega,L) = \inf_{C \cap \Omega
\neq \emptyset} \left\{\frac{\deg_L(C)}
{\sum_{x \in \Omega} \mult_x(C)}\right\};
$$
here the curve $C \subset X$   in the infimum is irreducible.  
We call a positive dimensional, irreducible
subvariety $V \subset X$ Seshadri exceptional for $L$
relative to $\Omega$ if
$$
\left(\frac{\deg_L(V)}{\sum_{x \in \Omega} \mult_x(V)}\right)^{\frac{1}
{\dim(V)}} = \epsilon(\Omega,L)
$$
and if $V$ is not properly contained in any subvariety $V^\prime$ satisfying
this same equality.  
\label{seshexc}
\end{definition}

\medskip

Let 
$G$ be  a connected
commutative algebraic group,  and $X$  a Serre compactification  \cite{S}  of $G$ (see \S \ref{subsecSerre}).  The group law on 
$G$ will play a fundamental role: it is used both to extend the collection of points
$\Omega = \{P_1, \ldots, P_n\}$ under consideration and to define the invariants which
enter into our main theorems.  
Using the language of Seshadri constants, our goal in this paper is to show that 
the interpolation problem for $X$, an ample line bundle $L$, and a finite subset $\Omega $ of $G$   can be reduced to
the corresponding interpolation problem on the compactifications of subgroups
of $G$ and their translates.  
More precisely, we will associate numerical invariants to each translate of a compactified
subgroup of $G$:
\begin{definition}  
Let 
$G$ be  a connected
commutative algebraic group,   $X$  a Serre compactification of $G$, 
  $L$  an ample line bundle on $X$,  and $\Omega  \subset G$ a finite subset.  Suppose $H  \subset G$ is a translate of a positive-dimensional connected subgroup and $\bar{H} \subset X$ 
its compactification.   Let
$$
\mu(\Omega,H,L) =  \frac{\left(
\frac{\deg_L(\bar{H})}{|(\Omega \cap H)(\dim(H))|}
\right)^{1/\dim(H)}}{\dim(H)}:
$$
here  $\Omega(\dim(H))$ denotes $\Omega + \ldots + \Omega$ with $\dim(H)$  summands.
Let
$$
\nu(\Omega,L) = \min_{H \subset G} \{\mu(\Omega, H ,L)\}
$$
where the minimum runs over all translates $H$ of  positive-dimensional connected subgroups   of
$G$.
\label{invariant}
\end{definition}

When Definition \ref{invariant} is applied in the special case where $H = G$
we will simply write
$$
\mu(\Omega,L) = \frac{\left(
\frac{L^d}{|\Omega(d)|}\right)^{1/d}}{d};
$$
throughout this paper we let $d=\dim (G)$.
Our first theorem shows that if $\epsilon(\Omega,L) < \mu(\Omega,L)$ then any Seshadri exceptional
subvariety $V$ for $\Omega$ and $L$ is contained in the compactification of a translate of a proper
algebraic subgroup of $G$, for which the corresponding inequality does not hold.  
\begin{theorem}  Suppose $L$ is an ample line bundle on $X$, a Serre
compactification of a connected commutative algebraic group $G$.  Suppose  $\Omega \subset G$ 
is a finite set 
such that
$$
\epsilon(\Omega,L) < \mu(\Omega,L).
$$
Let $V$ be a Seshadri exceptional subvariety of $X$ relative to $\Omega$, and let
$H$ denote the smallest translate of a connected
 algebraic subgroup  that contains $V\cap G$.
Then
\begin{itemize}
\item[$(i)$] The translate $H$ is distinct from $G$ and so  $V \cap G$ is degenerate.
\item[$(ii)$] We have $\epsilon(\Omega,L) \geq  \mu(\Omega,H,L) $.
\end{itemize}
\label{globaltheorem}
\end{theorem}

When the assumption $
\epsilon(\Omega,L) < \mu(\Omega,L) 
$ is not satisfied, part $(ii)$ of the conclusion holds with $H=G$. Therefore we obtain   as a corollary the following 
lower bound for the Seshadri constant $\epsilon(\Omega,L)$,
 which depends only on translates of
subgroups of $G$.

\begin{cor}
Suppose $L$ is ample on $X$, a Serre compactification of a connected
commutative algebraic
group $G$.  Suppose $\Omega \subset G$ 
is a finite set.
Then
$$
\epsilon(\Omega,L) \geq \nu(\Omega,L).
$$
\label{cor}
\end{cor}

Corollary \ref{cor} says that the {\em only} severe
obstructions for the interpolation problem along
 $\Omega$ come from translates of non--trivial
subgroups. This result is  close to being optimal, since linear algebra 
  (see  \cite{L} Proposition
5.1.9, stated below in Lemma \ref{except})  yields
$$
\epsilon(\Omega,L)  \leq  \min_{H \subset G}   \left(
\frac{\deg_L(\bar{H})}{| \Omega \cap H |}
\right)^{1/\dim(H)}  
$$
where the minimum runs over all translates $H$  of  positive-dimensional connected subgroups   of
$G$. Corollary  \ref{cor} is the best result one can hope for using the methods
of this paper.  It is as close to being optimal as the most
precise multiplicity estimate in this setting \cite{N3}.

\bigskip

A standard argument with vanishing theorems translates Corollary
\ref{cor}
into an interpolation estimate, but one which involves the canonical divisor $K_X$. Fortunately, 
this divisor is particularly convenient for the Serre compactifications of $G$: $-K_X$ is, 
up to algebraic equivalence,
 supported on $X \backslash G$ (see Lemma \ref{canonical} in \S \ref{subsecSerre} below) and effective,
 allowing us to remove the unwanted $K_X$:
\begin{cor}
Suppose $L$ is an ample line bundle
 on a Serre compactification $X$ of a connected commutative algebraic group 
$G$  of
dimension $d$.
Let $\Omega 
\subset G$ be a finite subset.  Suppose $\alpha$ is a positive integer
satisfying $d \leq \alpha < 
\nu(\Omega,L)$.  Then 
\begin{itemize}
\item[$(i)$]
The natural map
$$
H^0(X,K_X+L) \ra H^0\left(X, (K_X + L) \ts \oplus_{x \in \Omega}\oo_X/m_x^{\alpha+1-d}
 \right)
$$
is surjective.
\item[$(ii)$] The natural map
$$
H^0(X,L) \ra H^0\left(X,  L \ts \oplus_{x \in \Omega}\oo_X/m_x^{\alpha+1-d}
 \right)
$$
is surjective.
\end{itemize} 
\label{cor2}
\end{cor}
The hypothesis $ \nu(\Omega,L) > d$ will not in general be satisfied. 
But by \cite{L} Example 5.1.4
$\nu(\Omega,kL) = k\nu(\Omega,L)$  so
this can always be obtained by scaling $L$.
Note also that if $G$ is an abelian variety then $X=G$ is a Serre compactification, 
so that this corollary applies. 
These interpolation estimates are much more precise than those obtained by Masser \cite{M} and the first author \cite{Fi}; we refer to \S \ref{subseclien} for a detailed comparison. We will apply these new    estimates in a forthcoming paper \cite{FN}, in a situation where previously known results are not sufficiently precise.

\medskip

The basic method of proof employed in this article is that of \cite{NR}.
We would like to emphasize this similarity.  Multiplicity estimates and interpolation
estimates are established in the same way by studying the 
Seshadri exceptional subvariety
for $\Omega$ and $L$. The difference is that the obstruction subgroup for multiplicity
estimates can be larger than the obstruction for the interpolation problem. 
In \cite{N4} the
second author produced, for certain special cases of $\Omega$, a chain of subgroups,
the smallest of which is the obstruction to the interpolation problem and the 
largest of which is the obstruction to multiplicity estimates (see \cite{FN} for precise statements).

\medskip

The outline of the paper is as follows. We gather in \S \ref{secinterp} some properties of Serre compactifications and provide a detailed comparison with previously
known interpolation estimates. Then we move in \S \ref{secsesh} to   Seshadri constants and Seshadri
exceptional subvarieties, recalling their properties  and establishing how they behave under
translation in the case of interest here.  We also discuss differentiation of sections  of $L$
which is  fundamental to the proofs. Finally   in \S \ref{secpreuves}  we prove Theorem \ref{globaltheorem}
and derive Corollaries \ref{cor} and \ref{cor2}.

\medskip

\noindent
{\it Acknowledgments}  It is a pleasure to thank the Universit\'e de Paris-Sud, Orsay
for receiving the second author during January and February, 2012, providing an 
opportunity to start this work.  The second author would also like to thank Imperial
College which provided a pleasant environment in which to continue 
 work on this article. The first author is supported by Agence Nationale de la Recherche (project HAMOT), and would like to warmly
 thank Michel Waldschmidt for asking him
 to prove essentially the results presented here, almost 15 years ago in his thesis. The second author
would also like to thank Michel Waldschmidt for his long standing
 encouragement to work on the important and difficult
questions of multiplicity and interpolation estimates.

\section{Background} \label{secinterp}

In this section we recall the definition and properties of Serre  
compactifications (\S \ref{subsecSerre}) and study their canonical bundles. Then we compare in 
\S \ref{subseclien} our interpolation estimate to those of Masser and the first author; we mention there the case of arbitrary compactifications.

\subsection{Serre  Compactifications} \label{subsecSerre}

Except in \S \ref{subseclien} below, we shall always  compactify the commutative algebraic group $G$ following the
procedure due to   Serre \cite{S} that we describe now. Our proof might work for more general compactifications, but   technical difficulties would appear
that depend on the compactification. The same assumption appears already in \cite{N3}.
  
\bigskip
  
A commutative algebraic group $G$ can be viewed as an extension of an abelian variety by a linear
group.  In particular, there is an exact sequence of groups
$$
0 \ra L \ra G \ra A \ra 0
$$
where $A$ is an abelian variety and $L$ is a linear group.  The linear group $L$ can
be written as a product  $(G_m)^r \times (G_a)^s$  where $G_m$ is the multiplicative
group and $G_a$ the additive group: note that this expression of $L$ as a product is
{\em not} unique.  The linear group $L$ can be compactified as a product
of projective lines $(\pp^1)^{r+s}$.  The Serre compactification of $G$ is then
the induced compactifiation  $X$ obtained
by viewing $G$ as a principal fibre bundle over $A$ and then compactifying the fibres.  
To describe $X$ more concretely, suppose
$\pi: G \ra A$ is the projection map and $U \subset A$ is an open subset so that
$\pi^{-1}(U)$ can be written as a product $L \times U$.  If $p: X \ra A$ is the projection then
$p^{-1}(U)$ can be expressed as $(\pp^1)^{r+s} \times U$.   In other words, the linear
group $L$ can be compactified over an open cover of $A$ and then these open sets glue
together to give the Serre compactification $X$ of $G$.  It is important for us
that this compactification is equivariant, that is the group law on $G$ extends to
an action of $G$ on $X$.  Also important to note here is that $X \backslash G$
is a union of divisors.  The number of irreducible components of $X \backslash G$ is
$2r + s$ because $\pp^1 \backslash G_m$ is two points while $\pp^1 \backslash G_a$
is a single point.  

\bigskip

We shall compute now  the canonical divisor
$K_X$ on the Serre compactification of $G$ by reviewing 
  the argument at the end of \S 1 of \cite{N1}.  
In precise terms, we prove the following lemma which will allow  us to prove Corollary  \ref{cor2} at the end of \S \ref{secpreuves}.
  
\begin{lemma}
If $X$ is a Serre compactification of a commutative algebraic group $G$ then $-K_X$ is linearly equivalent to
the sum of a divisor $\pi^\ast(N)$, where $N$ is algebraically equivalent to zero on $A$, and
an effective divisor supported on $X \backslash G$. 
\label{canonical}
\end{lemma}

  \noindent
{\bf Proof of Lemma \ref{canonical}}     Knopf and Lange \cite{KL2} \S 2 show that if $M$ is an invertible sheaf 
on $X$ which admits an $L$ action  then $M$ can be expressed as
$$
M = \pi^\ast  (N) \ts  \left(\otimes_{i=1}^{r+s} \oo_X(a_iD_i)\right)
$$
where $N$ is an invertible sheaf  on $A$, the $a_i$ are integers,  
and the $D_i$ are the irreducible components of $X \backslash G$:  only
one of the two components coming from each compactifiation of $G_m$ is taken here.   Since the compactification 
 $X$ is equivariant, the canonical
sheaf $\oo_X(K_X)$ is acted on by $G$ hence by $L$.  Thus we can write
$$
\oo_X(K_X) = \pi^\ast(N)  \ts  \left(\otimes_{i=1}^{r+s} \oo_X(a_iD_i)\right).
$$
Since $G$ acts on $X$ and the canonical bundle of $X$ is preserved by this action,   it follows that $G$ acts on
and preserves $\pi^\ast (N)$.  This is only possible,
according to \cite{BL}  Corollary 2.5.4  or \cite{Mu} page 77 Theorem 1, if $N$ is algebraically equivalent to $0$.  
   To calculate
the integers $a_i$, we can consider $\oo_X(K_X)|p^{-1}(x)$ where $x$ is a point of $A$.  On the one
hand, $\oo_X(K_X)|p^{-1}(x) \simeq \oo_{p^{-1}(x)}(K_{p^{-1}(x)}$) as can be seen, for example, by
repeatedly applying the adjunction formula to divisors which are pulled back via $p$ from the base $A$.  On the
other hand, we know that $p^{-1}(x) = (\pp^1)^{r + s}$.  This tells us that $a_i = -2$ for all $i$, and   concludes the proof of Lemma \ref{canonical}.

\subsection{Connection with Previous Estimates} \label{subseclien}

In order to compare our results to that of \cite{M} and \cite{Fi}, 
  let us state and prove a weaker form a Corollary \ref{cor2}  valid for any    compactification of $G$. It involves a projective
embedding of $G$, which is   standard  when interpolation estimates   are applied to transcendental number theory.

\begin{cor}
Let $G$ be a connected commutative  algebraic group of dimension $d$ and $Y$   a  compactification of $G$. Let us fix a very ample divisor on $Y$, 
corresponding to   a locally closed immersion into a projective space
${\bf P}^N$. Then:
\begin{itemize}
\item[$(i)$]  There exists a positive constant $c$ (depending only on $Y$ and on this divisor) with the following property.
Let  $\Omega  
\subset G$ be a finite subset and $D, T$  positive integers such that
\begin{equation}
  D >c \,  T  \, |(\Omega\cap H)(\dim(H))|^{1/\dim(H)}
\label{inteqhypint}
\end{equation}
for any translate $H$ of a non-zero connected algebraic subgroup   of $G$. Then the evaluation map
$$
H^0(X,\oo(D)) \ra H^0\left(X, \oo(D) \ts \oplus_{x \in \Omega}\oo_X/m_x^{T}
 \right)
$$
is surjective.
\item[$(ii)$]  If $Y$ is a Serre compactification of $G$ then $(i)$ holds with $(\ref{inteqhypint})$ replaced with
\begin{equation}
(\deg_{\oo(1)} \bar H )  D ^{\dim(H)}>(\dim(H))^{\dim(H)}(T+d-1)^{\dim(H)}|(\Omega\cap H)(\dim(H))|.
\label{inteqhypintSerre}
\end{equation}
\end{itemize}
\label{cor4nv}
\end{cor}

\noindent
{\bf Proof of Corollary \ref{cor4nv}} Assertion $(ii)$ is just a reformulation of 
Corollary \ref{cor2} when $L$ is very ample. It implies that $(i)$ holds when $Y$ is a Serre compactification of $G$, with $c  = d^2+1$    in $(\ref{inteqhypint})$ (actually    any value greater than $d^2$
can be chosen, since $dT \geq T+d-1$).

Now let us fix a very ample divisor on  the Serre compactification of $G$, 
corresponding to   a locally closed immersion into  
${\bf P}^M$. In the setting of $(i)$,  consider  the identity map of $G$ to $G$  with respect to the embeddings of $G$ in  ${\bf P}^N$ on the left, and in  ${\bf P}^M$ on the right. On a open subset of $G$ containing $\Omega$, it is given by a family of $M+1$ homogeneous polynomials in $N+1$ variables, of the same degree $\delta$. Then $(i)$ holds in ${\bf P}^N$ with $c = \delta(d^2+1)$   
 in $(\ref{inteqhypint})$ because it does in ${\bf P}^M$ with $c = d^2+1$: this  follows from Lemma 1.3 and Proposition 2.3 of \cite{Fi} with $H=0$, $I' = I$, $D''= 0$ because conclusion $(ii)$ of this proposition never holds in this case. Note that in  \cite{Fi}   the set $\Omega$ is assumed to have a special form (namely $\Gamma(S_1,\ldots,S_r)$, see below), but this plays no role in this part of the paper  \cite{Fi}. This concludes the proof of Corollary \ref{cor4nv}.

\bigskip

We would like to emphasize the fact that part $(i)$ of Corollary \ref{cor4nv} is a drastic weakening of  Corollary \ref{cor2}, as the proof shows. We recall that $(ii)$,  as Corollary \ref{cor2}, is very precise: the best possible assumption  (up to
$O(D^{\dim(H)-1}) $) would be to replace $(\dim(H))^{\dim(H)} $ with 1 and $(\Omega\cap H)(\dim(H))$ with $\Omega\cap H$ in (\ref{inteqhypintSerre}).   It is essentially as precise as the second author's zero estimate  \cite{N3}. 

On the other hand, in $(i)$ the degree of $H$ does not appear in  (\ref{inteqhypint}), and 
  the constant $c$  depends  on the embedding of $G$ in ${\bf P}^N$ (as in early zero estimates, for instance \cite{MWun}); of course this is very unpleasant, and in contrast our results stated in the introduction depend only on $G$.
  In particular, unless $c$ is bounded explicitly, this corollary is only asymptotic with respect to $D$   :  (\ref{inteqhypint}) can hold only if $D$ is sufficently large with respect to $c$.

\bigskip

The previously known interpolation estimates in this setting (\cite{M}, \cite{Fi}) can be stated as follows: part $(i)$ of Corollary \ref{cor4nv} holds if 
  $\Omega = \Gamma(S_1,\ldots,S_r)$ is the set of elements $n_1\gamma_1+\ldots+n_r\gamma_r$ with fixed $\gamma_1,\ldots,\gamma_r\in G$ and integers $n_j$ with $|n_j| \leq  S_j$; moreover the constant $c$ in   (\ref{inteqhypint})   may  depend on  the embedding of $G$ in ${\bf P}^N$ and on $\gamma_1,\ldots,\gamma_r$. Therefore Corollary \ref{cor4nv}  contains, and refines, these results -- except  that in \cite{Fi} the order up to which derivatives are considered may vary according to the direction.
  
Whereas the proof of $(i)$ of Corollary \ref{cor4nv}  provides a value for $c$ (even though not a very natural one), the proof of \cite{M} and \cite{Fi} does not provide easily an explicit value. This constant 
    depends  (among others) on the degree of a family of homogeneous polynomials representing the addition law in this embedding (see  \cite{MWun}, pp. 492--494), and on the arbitrary choice of projective embeddings of finitely many quotients $G/H$ (see \cite{M}, p. 165, or \cite{Fi}, \S 2.1).
  
  A more serious drawback of the interpolation estimates of  \cite{M} and \cite{Fi} is that 
  $c$ depends also on $\gamma_1,\ldots,\gamma_r$. This makes the result  only asymptotic in $S_1,\ldots,S_r$ (and not only with respect to $D$, as $(i)$ of Corollary \ref{cor4nv}).
   For instance, if $\Omega$ has no special structure with respect to the group law of $G$, then one has to take $\{\gamma_1,\ldots,\gamma_r\} = \Omega$ and $S_1= \ldots = S_r = 1$ to apply it: the assumption  (\ref{inteqhypint})  means that $D/T$ is sufficiently large in terms of $\Omega$, and the result is completely trivial (it simply states that the Seshadri constant exists). The same conclusion holds if $\Omega$ consists only of torsion points (whereas this could be interesting in possible Diophantine applications); more generally, the torsion part of the $\Z$-module generated by $\Omega$ disappears in this result (see \cite{M}, \S 4, or \cite{Fi}, Step 2 in \S 4.2).

Finally,  Theorem \ref{globaltheorem} and Corollary \ref{cor}   provide a so-called {\em obstruction subgroup for interpolation} (in a terminology close to that of \cite{Fi}), that is a translate $H$ such that  $\epsilon(\Omega,L) \geq  \mu(\Omega,H,L) $, with an additional property: $H$ is the smallest translate that contains a Seshadri exceptional subvariety for $\Omega$ and $L$
(see also  \cite{NR} et \cite{N4}). On the contrary, no such interpretation follows from the proofs of  \cite{M} and \cite{Fi}.

\bigskip

There is only one  aspect of the first 
author's result \cite{Fi} which is   not contained in ours:   the fact that in \cite{Fi} the order up to which derivatives are considered may vary according to the direction. 
We would like to briefly address why the methods of this article do not suffice to establish an interpolation estimate in the
case of multiplicity along a proper
 analytic subgroup $0 \neq \Lambda \subset T_0(G)$.   In the multiplicity setting  where $\Lambda = T_0(G)$
 the jets which multiples of $L$ generate are controlled, up to requiring an additional $K_X$, by the Seshadri constant
 $\epsilon(\Omega,L)$ which has nice geometric properties.  When $\Lambda \neq T_0(G)$, there is no longer a simple global
 geometric invariant which measures positivity in the direction of $\Lambda$.  The Seshadri constant can be studied naturally
 on a single blow--up while the corresponding constant associated to derivation in the directions of $\Lambda$ requires more
 and more blow--ups as the order of jets increases.  These blow--ups then influence the end result adversely:  indeed, in Corollary \ref{cor2}, the $1-d$ in the exponent comes from the relative canonical bundle
of the blow up of $X$ along $\Omega$.    
In the case where $\Lambda \subset T_0(G)$
 is a proper subspace, the number of blow--ups is roughly equal to the order of vanishing on $\Lambda$ and each blow--up
 introduces a new exceptional divisor and a more complicated relative canonical bundle.   The fundamental issue here is that
 the canonical bundle $K_X$ which enters in these vanishing theorems is a global object associated to $X$, whereas what would
 be needed, in the case of a proper subspace $\Lambda \subset T_0(G)$, is a different object which only measures positivity
 in certain directions.   
   
\bigskip

Even though our point of view in this paper is to work only with Serre compactifications (because they make everything easier), it could be interesting 
to prove sharp interpolation estimates for other compactifications. Then  several parts
of our argument would have to change.  First,  we would need to
 define the translation operators $t_g$,  in a given embedding, in terms of homogeneous forms of bounded
degree: this unfortunately brings back all of the constants which depend on the embedding.  
Secondly and more seriously, we could no longer use
the fact that numerically equivalent line bundles have the same Seshadri exceptional subvarieties: we apply this to the line
bundle $L$ and its translates but in the absence of translation operators an alternative method would need to be developed
here.  Thirdly, an alternative means of differentiation would need to be found.

\bigskip

To conclude this section, we would  like to compare and contrast the method of proof here with the one employed
by Masser \cite{M} and generalized by the first author \cite{Fi}.
Rather than using positivity and
vanishing theorems as is done here, Masser  takes a more concrete approach we summarize briefly here.
If the evaluation map is not surjective (say with no multiplicities, i.e. $T=1$) and if there is no obstructing subgroup,  then there is a linear relation between the values at the points of $\Omega$ 
of all $P \in H^0(X,\oo(D))$ -- that is, a non-zero functional that vanishes when applied to any $P$. Then translating this functional yields many relations between the values at the points of $\Omega(d)$. Linear algebra yields 
the existence of a  non-zero $P \in H^0(X,\oo(D'))$ (for some $D'$ close to $D$) which vanishes at  many points among those of $\Omega(d)$; these many relations imply that it vanishes at all points of $\Omega(d)$. This contradicts a zero estimate, for instance that of Philippon~\cite{P1}.

\section{Seshadri Constants and Seshadri Exceptional Subvarieties} \label{secsesh}

In this section we gather together all of the preliminary results which we will need
in proving the main theorems. 
We first recall in \S \ref{subsecsesh1}  some of the important properties of Seshadri 
constants and Seshadri 
exceptional subvarieties which we use repeatedly in our proofs. Then we focus in  \S \ref{subsecsesh2} on the order of vanishing of sections along 
Seshadri 
exceptional subvarieties: here we prove Theorem \ref{vanishing}, the main ingredient in our approach.

\subsection{General Properties} \label{subsecsesh1}

Most  results of this section apply in general  and do not require that $X$ be a  compactification 
of a commutative algebraic group. 
Until the end of this section, we let $X$ denote a smooth projective variety with additional hypotheses as needed. These results will be applied in \S \ref{secpreuves} to the case where
$X$ is a Serre compactification of a commutative algebraic group, but they hold in general and are of independent interest.

\medskip

We first discuss in detail Seshadri constants and Seshadri exceptional subvarieties.
A useful reference for this material is Chapter 5 of \cite{L}, 
where complete proofs are provided when $\Omega$
consists in a single point: the general results for a finite subset $\Omega$ 
can be proved in the same way.  In addition to the
definition of Seshadri constants in terms of curves and multiplicities which was given in the introduction,
there is an alternative definition which is often useful.  
For this definition, given in \cite{L} Definition 5.1.1, note that a divisor
 $D$ on $X$ is called
{\it nef} if $D \cdot C \geq 0$ for every curve $C \subset X$.
\begin{definition}
Suppose $\Omega \subset X$ is a finite set, $L$ an ample line bundle on $X$.  Let $\pi: Y \ra X$ be
the blow--up of $X$ along $\Omega$ with exceptional divisor $E$.  Then the Seshadri constant of 
$L$ along $\Omega$ is$$
\epsilon(\Omega,L) = \sup\{\alpha \geq 0: \pi^\ast(L)(-\alpha E) \,\,\,\mbox{is nef}\}.
$$
\label{sesh2}
\end{definition}
The fact that Definition \ref{seshexc} and Definition \ref{sesh2} are equivalent can be found in \cite{L} Proposition 5.1.5. 
We will only use this alternative characterization of Seshadri constants in the proof of Corollary 
\ref{cor2}.
Campana and Peternell \cite{CP} established that a Seshadri exceptional subvariety for $L$ and $\Omega$, as defined in Definition \ref{seshexc}, 
always exists;  a proof can be found in \cite{L} Proposition 2.3.18.  
Combining  Definition \ref{seshexc} with Proposition 5.1.9 of \cite{L} yields the following lemma:
\begin{lemma}
For any positive-dimensional irreducible subvariety $V$ of $X$ we have
$$
\epsilon(\Omega,L)^{\dim(V)}\left(\sum_{x \in \Omega \cap V} \mult_x(V)\right) \leq \deg_L(V).
$$
If  $V$ is a Seshadri exceptional subvariety with respect to $L$ and $\Omega$, then equality holds.  
\label{except}
\end{lemma}

Of great
importance to us is the fact that Seshadri exceptional subvarieties  for $\Omega$ and $L$ depend  only on the
numerical equivalence class of $L$: recall that two line bundles $L$ and $M$ on a variety
$X$ are called numerically equivalent if $L \cdot C = M \cdot C$ for every curve
$C \subset X$.  
\begin{lemma}  Suppose  $M$ and $L$ are numerically equivalent ample  line bundles on a smooth projective variety $X$. Let
$\Omega \subset X$ be a finite subset.   Then a  subvariety $V \subset X$ is a Seshadri exceptional for $\Omega$ and $L$
if and only if it is Seshadri exceptional for $\Omega$ and $M$.
\label{num}
\end{lemma}

\noindent
{\bf Proof of Lemma \ref{num}} Using Definition \ref{seshexc} 
of $\epsilon(\Omega,L)$ it is clear that $\epsilon(\Omega,L) = \epsilon(\Omega,M)$ when
$L$ and $M$ are numerically equivalent because 
$$
\frac{L \cdot C}{\sum_{x \in \Omega} \mult_x C} = \frac{M \cdot C}{\sum_{x \in \Omega} \mult_x C}
$$
for every curve $C$ which contains at least one point of $\Omega$.  
If $V$ is a subvariety in $X$ which is Seshadri
exceptional with respect to $L$ and $\Omega$, Definition \ref{seshexc} 
shows  that $V$ is also  Seshadri exceptional relative
to $M$ and $\Omega$ (and vice versa) as desired: here we use the fact that degrees of all subvarieties, not
just curves, are equal with respect to numerically equivalent line bundles.

\medskip

Lemma \ref{num} has an important application in our situation.  Because
$X$ is an equivariant compactification of $G$ there
is a morphism $f: G \times X \rightarrow X$ so that, restricted to $G \times G$,
$f$ gives the addition law on $G$. 
If $p_1: G \times X  \ra G$ denotes the projection
to the first factor then, by \cite{Fu} Definition 10.3, the line bundles
$f^\ast(L)|p_1^{-1}(g)$ are algebraically equivalent for all $g \in G$.  
Each fibre $p_1^{-1}(g)$ is isomorphic to $X$ and, via this isomorphism,
$f: p_1^{-1}(g) \ra X$ is identified with $t_g: X \ra X$, the translation map
given by $t_g(x) = g + x$ for all $x \in X$.
We conclude that  the line bundles $\{t_g^\ast(L)\}_{g \in G}$ on $X$ are all algebraically equivalent
to one another. 
By \cite{Fu} \S 19.1 we deduce that the line bundles
$t_g^\ast(L)$ are all numerically equivalent.  Hence, we derive from Lemma \ref{num} the
important corollary:
\begin{cor}
For any  ample  line bundle $L$ on $X$ and any $g \in G$, a subvariety $V$ is Seshadri exceptional for $\Omega$ and $L$ if and
only if it is Seshadri exceptional for $\Omega$ and  $t_g^\ast(L)$.
\label{translates}
\end{cor}

\medskip

Before presenting in \S \ref{subsecsesh2} the main result of this section, we require one more preliminary lemma about jet separation.
Using the notation of \cite{L} Definition 5.1.15,
 write $J_x^s(L)$ for $H^0(X,L \ts \oo_X/m_x^{s+1})$, the jets of
order $s$ at $x$ for a line bundle $L$.  Similarly, 
 write $J_\Omega^s(L) = \oplus_{x \in \Omega} J_x^s(L)$.  We say that a line bundle $L$ separates $s$--jets
 along $\Omega$ if the natural map
 $$
 H^0(X,L) \ra \bigoplus_{x \in \Omega} J_x^s(L) 
 $$
 is surjective.
Following \cite{L} Definition 5.1.16 we
write $s(L,\Omega)$ for the largest non--negative integer $s$  such that $L$ separates 
$s$--jets along $\Omega$, assigning
the value $-1$ if $L$ does not separate zero jets along $\Omega$.
 Proposition
5.1.17 of  \cite{L} gives the following important relationship between the Seshadri constant $\epsilon(\Omega,L)$ and
the asymptotic separation of jets by powers of the line bundle $L$:
\begin{lemma}
Suppose $L$ is an ample line bundle on a smooth projective variety $X$ and $\Omega \subset X$ a finite subset.  Then
$$
\epsilon(\Omega,L) = \lim_{k \ra \infty}\frac{s(kL,\Omega)}{k}.
$$
\label{separation}
\end{lemma}

\subsection{Order of Vanishing along Exceptional Subvarieties} \label{subsecsesh2}

From now on we denote by $X$ a Serre  compactification of a connected commutative  algebraic group $G$.
To begin with, let us state (and prove)  the main theorem of this section, relating the order of vanishing of a section of $L$ along
$\Omega$ to its order of vanishing along any Seshadri 
exceptional subvariety $V$ relative to $L$ and $\Omega$.
It  is closely related to \cite{ekl} Proposition 2.3, \cite{N5} Lemma 1.3, 
and is
a simplified version of \cite{NR} Proposition 5:

\begin{theorem}
Suppose $X$ is a Serre compactification of a connected commutative  algebraic group $G$.  
Let   $L$ be
an ample line bundle on $X$,  $\Omega \subset X$ a finite subset, $s \in H^0(X,L)$ a non--zero
section.  Suppose $V$ is Seshadri exceptional for $L$ relative to $\Omega$.  Let 
$m = \mult_\Omega(s)$.  Then
$$
\mult_V(s) \geq m - \epsilon(\Omega,L).
$$
\label{vanishing}
\end{theorem}

An important feature in our proof is that we differentiate sections, in the same way as in   \cite{N3,NR}; let us recall this briefly.
To each
non-zero vector $v \in T_0(G)$ is associated a translation invariant vector field on $X$. 
If $s \in H^0(X, L)$
then the derivative of $s$ can be taken locally
with respect to this vector field and we denote this derivative by $D_v(s)$.  
This is not in general a section of $L$ on $X$ because the local patching restraints on $s$ are destroyed
when taking derivatives:  however, $D_v(s)$ is a well--defined section of $H^0(Z,L)$ if $Z$ is contained
in the zero locus of $s$.  For any $r \geq 1$ we shall consider differential operators of order $r$, that is polynomials of degree $r$ in these operators $D_v$. In the same way,  if $s \in H^0(X, L)$
 vanishes along $Z$ with multiplicity at least $\mu$ then $D(s) \in H^0(Z,L)$ for any differential operator $D$ of order less than or equal to $\mu$.  The reader is referred
to \S 1 of \cite{N3} for more details.
We shall use the following result, which is a special case of  Lemma 4 from \cite{N3}.

\begin{lemma}
Suppose $X$ is a Serre compactification of a connected commutative  algebraic group $G$.
  Let $L$ be an ample line bundle on $X$.
Let $x\in X$ and  $0 \neq s \in H^0(X,L)$. Let 
$V \subset X$ be an irreducible subvariety  containing $x$, and   $D$ be a differential operator of
order $r $ on $X$. Assume that $r =  \mult_V(s) < \mult_x(s) = M$.
 Let $Z(D(s))  = \sum_{i=1}^n a_iW_i$ where $D(s)$ is considered as a global section of $L$ on $V$.
Then
$$
\sum_{i=1}^n a_i\mult_x(W_i) \geq \mult_x(V)(M-r).
$$
\label{derivation}
\end{lemma}

\noindent {\bf Proof of Theorem \ref{vanishing}}  Suppose that the conclusion of Theorem \ref{vanishing} is false so
that
$$
\mult_V(s) < m - \epsilon(\Omega,L).
$$
Choose a differential operator $D$ on $X$
  of order  $ \mult_V(s) $   so that $D(s) \in H^0(V,L)$ is well defined and non--zero.  
Write
$$
Z(D(s)) =\sum_{i=1}^n a_iW_i.
$$
Lemma \ref{derivation} gives
\rn
\begin{equation}
\mult_x(D(s)) =  \sum_{i=1}^n a_i \mult_x(W_i) > \mult_x(V)\epsilon(\Omega,L) \,\,\,\mbox{for each $x \in \Omega \cap V$. }
\label{plop}
\end{equation}
 Applying Lemma \ref{separation} and (\ref{plop}), we may  choose $k$ sufficiently large so that for each $x \in \Omega \cap V$  \rn
 \begin{equation}
 \mult_x(D(s))\left(\frac{s(kL,\Omega)}{k}\right)^{\dim (V)-1} > \mult_x(V) \epsilon(\Omega,L)^{\dim (V)}.
 \label{dope}
 \end{equation}
 By definition of $s(kL,\Omega)$ there exist sections $s_1, \ldots, s_{\dim(V) - 1} \in H^0(X,kL)$
 so that 
\rn
\begin{equation}
\mult_x(s_i) = s(kL,\Omega)
\label{hope}
\end{equation}
 for each $x \in \Omega$ and the tangent cones (see \cite{Fu} page 227) of the divisors $Z(s_i)$ meet properly  at
 each  point in $x \in \Omega \cap V$.  It follows that each $x \in \Omega \cap V$  is an isolated
irreducible component of 
\rn
\begin{equation}
Z(s_1) \cap \ldots \cap Z(s_{\dim(V) -1}) \cap Z(D(s)). 
\label{elope}
\end{equation}

Although each point of $x \in \Omega \cap V$ is an irreducible component of
the intersection (\ref{elope}), there may also be positive dimensional
components $\{W'_j\}$ in (\ref{elope}) which do not contain any point of $\Omega \cap V$.  
Since $L$ is an ample line bundle, $\deg_{kL}(W'_j) > 0$ for all $j$ and
so the part of the intersection (\ref{elope}) which is not supported on 
$\Omega \cap V$ can be represented by an effective zero cycle. 
By \cite{Fu} Corollary 12.4, for each $x \in \Omega \cap V$
the multiplicity of $x$ in $Z(s_1) \cap \ldots \cap Z(s_{\dim(V)-1}) \cap Z(D(s))$ is at least
$$
\mult_xZ(D(s)) \left( \prod_{i=1}^{\dim(V) -1} \mult_x(Z(s_i)) \right).
$$
 Hence we find
 \begin{eqnarray*}
 k^{\dim(V)-1}\deg_L(V) &=& \deg \left(Z(s_1) \cap \ldots \cap Z(s_{\dim(V)-1}) \cap Z(D(s))\right) \\
                               &\geq&  \sum_{x \in \Omega \cap V}  \mult_xZ(D(s)) \left( \prod_{i=1}^{\dim(V) -1} \mult_x(Z(s_i)) \right) \\
                                 &>&  \left(\sum_{x \in \Omega \cap V} \mult_x(V)\right)k^{\dim(V)-1}\epsilon(\Omega,L)^{\dim(V)}
 \end{eqnarray*}
 with the last inequality coming from (\ref{dope}) and (\ref{hope}).
This contradicts Lemma \ref{except} and establishes Theorem \ref{vanishing}. \
\bigskip

ln order to use Theorem \ref{vanishing} in a broader setting in what follows we require 
\begin{cor}
Suppose $X$ is a Serre compactification of a connected commutative  algebraic group $G$.  Suppose  
 $L$ an ample line bundle on $X$,  $\Omega \subset X$ a finite subset, $s \in H^0(X,L)$ a non--zero
section.  Suppose $V$ is Seshadri exceptional for $L$ relative to $\Omega$.  Let  $g \in G$ and 
$m = \mult_{g+\Omega}(s)$.  Then
$$
\mult_{g+V}(s) \geq m - \epsilon(\Omega,L).
$$
\label{vanishing1}
\end{cor}

\noindent
{\bf Proof of Corollary \ref{vanishing1}}  Let $t_{g}: X \ra X$ denote the map given by 
$t_{g}(x) = g + x$ for all $x \in X$.  It follows from Corollary \ref{translates}
that $V$ is Seshadri exceptional for $t_{g}^\ast(L)$
relative to $\Omega$.  If $s$ is the section of Corollary \ref{vanishing1} then
$t_g^\ast(s)$ vanishes to order  $m$ along $\Omega$.  Since $V$ is Seshadri exceptional
for $\Omega$ and $t_g^\ast(L)$ we can apply Theorem \ref{vanishing} to conclude that
$$
\mult_V(t_{g}^\ast(s)) \geq m - \epsilon(\Omega,L).
$$
But $\mult_V(t_g^\ast(s)) = \mult_{g + V}(s)$ so this concludes the proof of Corollary \ref{vanishing1}.

\section{Proofs of Main Results} \label{secpreuves}

In this section we prove Theorem \ref{globaltheorem} and its two Corollaries \ref{cor} and \ref{cor2}.

\bigskip

\noindent
{\bf Proof of Theorem \ref{globaltheorem}}  
By the asymptotic Riemann--Roch Theorem, \cite{L} Example 1.2.19,
$$
h^0(X,kL) = \frac{L^dk^d}{d!} + O(k^{d-1}).
$$
If $\alpha$ is a real number, we let $\lceil k\alpha\rceil$ denote the smallest integer which is greater than or equal to
$k\alpha$.  We have
$$
\dim\left(J_{\Omega(d)}^{\lceil k\alpha\rceil}\right) = \frac{(k \alpha)^d}{d!}|\Omega(d)|  + O(k^{d-1}).
$$
Hence as long as
$|\Omega(d)|\alpha^d < L^d$ or
\rn
\begin{equation}
\alpha < \left(\frac{L^d}{|\Omega(d)|}\right)^{\frac{1}{d}}, 
\label{cope}
\end{equation}
when $k \gg 0$ there is a non--zero section $s \in H^0(X,kL)$ 
whose jets of order $\lceil k\alpha \rceil$ are identically zero. Consequently,  
 there exists a non--zero section $s \in H^0(X,kL)$ such that  
\rn
\begin{equation}
\mult_{\Omega(d)}(s) > k\alpha
\label{room}.
\end{equation}
Since 
$$
\epsilon(\Omega,L) < 
\mu(\Omega,L) = \frac{\left(\frac{L^d}{|\Omega(d)|}\right)^{\frac{1}{d}}}{d}
$$ 
we can apply (\ref{cope})
and (\ref{room}) to conclude that for $k$ sufficiently large
there is a non--zero section $s \in H^0(X,kL)$ such that 
\rn
\begin{equation}
\mult_{\Omega(d)}(s) > kd\epsilon(\Omega,L).
\label{boom}
\end{equation}

Since $V$ is Seshadri exceptional for $L$ and $\Omega$, $V$ must contain at
least one point of $\Omega$ and so $V \cap G$ is non--empty.  
 Let $V^{( r)} = (V \cap G)  + \ldots + (V \cap G)$
with $r$ summands.  Let $V^{(0)} = \{0\}$ and $\Omega(0) = \{0\}$ where $0 \in G$ is the identity
element.  
We will show by induction on $r$ that for $0 \leq r \leq d$ we have, provided  $k$ is  sufficiently large:
\rn
 \begin{equation}
\mult_{\Omega(d-r) + V^{( r)}}(s) > k(d-r)\epsilon(\Omega,L).
\label{chain}
\end{equation}
For $r = 0$, (\ref{chain}) is (\ref{boom}).
Suppose $ 0 \leq r \leq d-1$ and that
(\ref{chain}) has been verified for $r$.  
For each $x \in \Omega(d-r-1) + V^{( r)}$,
we have $x + \Omega \subset \Omega(d-r) + V^{( r)}$.  By (\ref{chain})
$$
\mult_{x+ \Omega}(s) > k(d-r)\epsilon(\Omega,L)  \,\,\,\mbox{ for any }  x \in \Omega(d-r-1) + V^{( r)}.
$$
By Corollary \ref{vanishing1}
$$
\mult_{x  + V}(s) > k(d-r-1)\epsilon(\Omega,L)  \,\,\, \mbox{ for any } x \in \Omega(d-r-1) + V^{( r)}.
$$
Thus
$$
\mult_{x + V^{(r+1)}}(s) > k(d - r - 1)\epsilon(\Omega,L)   \,\,\, \mbox{ for any } x \in \Omega(d-r-1)
$$
and this is exactly (\ref{chain}) for the case $r+1$.
When  $r = d$, we conclude
that $s$ vanishes along $(V \cap G) + \ldots + (V \cap G)$ with $d$ summands.  Since
$s$ is non--zero and $d = \dim(G)$ this is only possible if $V \cap G$ is 
contained in a  translate of a proper connected algebraic subgroup of $G$. This concludes the proof of  $(i)$ in Theorem \ref{globaltheorem}.

\medskip

Suppose that $(ii)$ is false so that
$$
\epsilon(\Omega,L) < \mu(\Omega,H,L).
$$
We will write $Y = \bar{H}$, the Zariski closure of $H$ in $X$.
We claim that $V$ is a Seshadri exceptional subvariety relative to
 $L|Y$ and $\Omega \cap H$, 
where $L|Y$ denotes the restriction of $L$ to $Y$.   
To see this note that $H = Y \cap G$ is a 
translate of a connected algebraic subgroup and so is smooth. Consequently
 it makes sense to talk about the Seshadri constant $\epsilon(\Omega \cap Y,L|Y)$. 
Thinking of $\epsilon(\Omega \cap Y,L|Y)$ in terms of Definition \ref{seshexc}, we have
\rn
\begin{equation}
 \epsilon(\Omega \cap Y,L|Y) \geq \epsilon(\Omega, L)
 \label{limp}
\end{equation}
 because the infimum on the left
is taken over a smaller collection of curves than the infimum on the right,
and the multiplicity of these curves is added on
 a smaller set of points.  
 On the other hand, using Definition \ref{seshexc} we have
 $$
\epsilon(\Omega,L) =  \left(\frac{\deg_L(V)}{\sum_{x \in \Omega} \mult_x(V)}\right)^{\frac{1}{\dim(V)}}.
 $$
 But $V\cap G \subset H$ and consequently the only points of $\Omega$ which can be in $V$ are those that are also
 in $H$ so we find
 $$
 \epsilon(\Omega,L) =  \left(\frac{\deg_L(V)}{\sum_{x \in \Omega \cap H} \mult_x(V)}\right)^{\frac{1}{\dim(V)}}.
 $$
 Now $V = \overline{V\cap G}$  so that
 Lemma \ref{except} yields 
 $$
 \epsilon(\Omega \cap Y, L|Y) \leq \left(\frac{\deg_L(V)}{\sum_{x \in \Omega \cap H} \mult_x(V)}\right)^{\frac{1}{\dim(V)}} = \epsilon(\Omega,L)  
$$
 and thus, recalling (\ref{limp}),  equality must hold here and in (\ref{limp}).  
Hence $V$ is Seshadri exceptional for $L|Y$ and $\Omega \cap H$.

 We have assumed that $\epsilon(\Omega,L) < \mu(\Omega,H,L)$.  Since $Y=\bar H$ and  equality   holds   in (\ref{limp})  this means  
\begin{equation} \label{equtile}
 \epsilon(\Omega \cap H,L|Y) <  \frac{\left(
\frac{\deg_L(Y)}{|(\Omega \cap H)(\dim(H))|}
\right)^{1/\dim(H)}}{\dim(H)}.
\end{equation}
 If $H \subset G$ is a subgroup, then we may apply part $(i)$ to $H$, $L|Y$, and $\Omega \cap H$
 and this leads to a contradiction. 
Suppose then that $H$ is a translate of a connected, proper algebraic subgroup $H_0 \subset G$. 
Choose $g \in \Omega \cap H$ and
 let $t_g: X \ra X$ denote
 translation by $g$.  
 Denote by 
 $$
 \tau_g: \overline{H_0} \ra \overline{H}
 $$
 the restriction of $t_g$ to $\overline{H_0}$.
   
 Since $t_g$ is an isomorphism it preserves Sesahdri exceptional subvarieties:
 $t_{-g}(V)$ is Seshadri exceptional relative to $t_{-g}(\Omega)$ and $t_g^\ast(L)$ and
 $\epsilon(t_{-g}(\Omega),t_g^\ast(L)) = \epsilon(\Omega,L)$.  By Corollary
 \ref{translates}, $t_{-g}(V)$ is Seshadri exceptional relative to $t_{-g}(\Omega)$ and $L$.  
 In the same way,  $\tau_{-g}(V)$ is Seshadri
 exceptional relative to $\tau_{-g}(\Omega \cap H) = t_{-g}(\Omega)\cap H_0$ and $L|\overline{H_0}$
 because  $V$ is Seshadri exceptional for $L|Y$ and $\Omega \cap H$, and   we have
 \begin{equation}
 \epsilon(t_{-g}(\Omega) \cap H_0, L|\overline{H_0}) = \epsilon(\Omega \cap H,L|\overline{H}).
 \label{1}
 \end{equation}
Looking at Defnition \ref{seshexc}, since $\deg_L(\overline{H_0}) = \deg_L(\overline{H})$ and $|(t_{-g}(\Omega) \cap H_0)
(\dim(H_0)| =    |(\Omega \cap H) (\dim(H))|$ we find
 \begin{equation}
 \mu(\Omega,H,L) =\mu(t_{-g}(\Omega)\cap H_0, H_0,L).
 \label{2}
 \end{equation}
Combining (\ref{equtile}), (\ref{1}), and (\ref{2}),   part $(i)$ of Theorem \ref{globaltheorem} shows that $\tau_{-g}(V)$ is degenerate in $\overline{H_0}$ which
is not possible since $H_0$ is, by hypothesis, the smallest connected subgroup of $G$ containing $t_{-g}(V)$.   This concludes
the proof of Theorem \ref{globaltheorem}.
 
 \bigskip 
 
 \noindent
 {\bf Proof of Corollary \ref{cor}}
 By Theorem \ref{globaltheorem}, either $\epsilon(\Omega,L) \geq \mu(\Omega,L)$ or  there is a translate $H$  of a connected 
 proper algebraic subgroup   so that $\epsilon(\Omega,L) \geq \mu(\Omega,H,L)$.  Since $\mu(\Omega,L) \geq \nu(\Omega,L)$
 and $\mu(\Omega,H,L) \geq \nu(\Omega,L)$  by hypothesis, we conclude that
 $$
 \epsilon(\Omega,L) \geq \nu(\Omega,L)
 $$
 as desired.
 
 \bigskip
 
  \noindent
 {\bf Proof of Corollary \ref{cor2}}
Let $\pi: Y \ra X$ be the blow up of
$\Omega$ with exceptional divisor $E$.
By Definition \ref{sesh2},
$\pi^\ast(L)\left(-\epsilon(\Omega,L) E\right)$  is a limit of nef line bundles and 
so \cite{L} Example 1.4.16 shows that
$\pi^\ast(L)\left(-\epsilon(\Omega,L) E\right)$ is nef.  Thus $$
 \pi^\ast(L)\left(-\epsilon(\Omega,L) E\right)^d \geq 0.
$$
On the other hand, by hypothesis,  $0 < \alpha < \nu(\Omega,L)$ so by Corollary \ref{cor} we have $\alpha < \epsilon(\Omega,L)$.  Thus
$$
\pi^\ast(L)(-\alpha E)^{d}  >  \pi^\ast(L)(-\epsilon(\Omega,L) E)^d \geq 0.
$$
 By  \cite{L} Theorem 2.2.16, $\pi^\ast(L)(-\alpha E)$ is a  big line
bundle (see \cite{L} Definition 2.2.1). 
We may therefore apply the Kawamata--Viehweg vanishing 
theorem, \cite{L} Theorem 4.3.1, and obtain
$$
H^1\left(Y,(K_Y+\pi^\ast(L))\left(-\alpha E\right)\right) = 0.
$$
Since $K_Y =  \pi^\ast(K_X) + (d-1)E$ this means
\rn
\begin{equation}
H^1\left(Y,\pi^\ast(K_X + L)
\left((d-1-\alpha) E\right)\right) = 0.
\label{projection}
\end{equation}

Let $\ii_\Omega(\alpha +1-d) \subset \oo_X$ denote the sheaf of functions vanishing to
order at least $\alpha+2-d$ at each point of $\Omega$. 
Applying the projection formula (\cite{H} Exercise 8.1 and \cite{KMM} Theorem 1-2-3)
to (\ref{projection}) gives  \
\rn
\begin{equation}
H^1\left(X,K_X + L \ts \ii_\Omega(\alpha +1-d)\right) = 0.
\label{grap}
\end{equation}
Using (\ref{grap}) and the long exact cohomology sequence associated to the exact
sequence of sheaves
$$
0 \ra H^0(X, (K_X + L) \ts \ii_{\Omega}(s)) \ra H^0(X, K_X + L) \ra
J_\Omega^{\alpha+1-d}(K_X +L) \ra 0
$$
shows that the map
$$
H^0(X,K_X + L) \ra J_\Omega^{\alpha+1-d}(K_X +L)
$$
is surjective, concluding the proof of $(i)$ of Corollary \ref{cor2}.  

\bigskip

In the proof of Corollary \ref{cor2} $(i)$, 
the line bundle $\pi^\ast(L)(-\alpha E)$ is actually ample as is shown
in \cite{NR} page 923.  
As a result, the Kawamata--Viehweg vanishing theorem is not really needed and
the Kodaira Vanishing Theorem, \cite{L} Theorem 4.2.1, is sufficient.  
The important fact in the proof of Corollary \ref{cor2} $(i)$ is that $\pi^\ast(L)(-\epsilon(\Omega,L)E)$ is big and nef.
If $M$ is a line bundle on $X$ which is numerically equivalent to zero then $\pi^\ast(L+M)(-\epsilon(\Omega,L)E)$
is also big and nef and the proof of Corollary \ref{cor2} $(i)$  gives the surjection
 \rn
\begin{eqnarray}
H^0\left(X, K_X + L + M\right) \ra J_\Omega^{\alpha + 1 - d}(K_X +  L + M)
\label{grapde}
\end{eqnarray}
  which will be vital in the proof of Corollary \ref{cor2} $(ii)$ which we provide now.

\bigskip

 According to Lemma \ref{canonical}, $-K_X$ is linearly equivalent to the sum of an effective divisor $E$ supported 
on $X  \backslash G$ and $\pi^\ast(N)$ where $N$ is algebraically equivalent to zero.  Since $\pi^\ast(N)$ is
numerically equivalent to zero, (\ref{grapde}) with $M = \pi^\ast(N)$ gives a surjection
$$
H^0\left(X,  K_X+ L  + \pi^\ast(N)\right)  \ra J_\Omega^{\alpha + 1 - d}(K_X + L +\pi^\ast (N)).
$$
Using the fact that $K_X$ is linearly equivalent to  $- E -\pi^\ast(N)$  we obtain a surjection
$$
H^0(X, -E +L)  \ra J_\Omega^{\alpha + 1 - d}(- E+L).
$$
Since $E$ is supported on $X \backslash G$ and $\Omega \subset G$ 
adding the divisor $E$  induces the desired surjection
$$
H^0(X, L)  \ra J_\Omega^{\alpha + 1 - d}(L).
$$

   \medskip

\medskip
\medskip

\begin{tabular}{l}
St\'ephane Fischler\\
Univ Paris-Sud, Laboratoire de Math\'ematiques d'Orsay\\
CNRS,  F-91405 Orsay, France\\ 
{\em Electronic mail:} stephane.fischler@math.u-psud.fr
\end{tabular}

\medskip

\begin{tabular}{l}
Michael Nakamaye\\
Department of Mathematics and Statistics\\
University of New Mexico\\
Albuquerque, New Mexico 87131, U.S.A.\\
{\em Electronic mail:} nakamaye@math.unm.edu
\end{tabular}

\newcommand{\url}{\texttt}

\end{document}